\documentclass[12pt,twoside,leqno]{article}
\usepackage{amsmath}
\usepackage{amssymb}
\usepackage{amsxtra}
\usepackage{amscd}
\usepackage{amsthm}
\usepackage[mathscr]{eucal}

\theoremstyle{plain}
\newtheorem{thm}[subsection]{Theorem}

\newtheorem{sbthm}[subsubsection]{Theorem}
\newtheorem{sbprop}[subsubsection]{Proposition}

\newtheorem{sblem}[subsubsection]{Lemma}

\theoremstyle{definition}

\newtheorem{para}[subsection]{}

\newtheorem{sbrem}[subsubsection]{Remark}

\newtheorem{sbpara}[subsubsection]{}

\newenvironment{pf}{\proof[\proofname]}{\endproof}

\begin{document}

\title{On SL(2)-orbit theorems}

\author{Kazuya Kato}

\maketitle

\newcommand\Cal{\mathcal}
\newcommand\define{\newcommand}

\define\gp{\mathrm{gp}}%
\define\fs{\mathrm{fs}}%
\define\an{\mathrm{an}}%
\define\mult{\mathrm{mult}}%
\define\Ker{\mathrm{Ker}\,}%
\define\Coker{\mathrm{Coker}\,}%
\define\Hom{\mathrm{Hom}\,}%
\define\Ext{\mathrm{Ext}\,}%
\define\rank{\mathrm{rank}\,}%
\define\gr{\mathrm{gr}}%
\define\cHom{\Cal{Hom}}
\define\cExt{\Cal Ext\,}%

\define\cC{\Cal C}
\define\cD{\Cal D}
\define\cO{\Cal O}
\define\cS{\Cal S}
\define\cM{\Cal M}
\define\cG{\Cal G}
\define\cH{\Cal H}
\define\cE{\Cal E}
\define\cF{\Cal F}
\define\cN{\Cal N}
\define\fF{\frak F}
\define\Dc{\check{D}}
\define\Ec{\check{E}}

\newcommand{\N}{{\mathbb{N}}}
\newcommand{\Q}{{\mathbb{Q}}}
\newcommand{\Z}{{\mathbb{Z}}}
\newcommand{\R}{{\mathbb{R}}}
\newcommand{\C}{{\mathbb{C}}}
\newcommand{\bN}{{\mathbb{N}}}
\newcommand{\bQ}{{\mathbb{Q}}}
\newcommand{\bF}{{\mathbb{F}}}
\newcommand{\bZ}{{\mathbb{Z}}}
\newcommand{\bP}{{\mathbb{P}}}
\newcommand{\bR}{{\mathbb{R}}}
\newcommand{\bC}{{\mathbb{C}}}
\newcommand{\bbQ}{{\bar \mathbb{Q}}}
\newcommand{\ol}[1]{\overline{#1}}
\newcommand{\too}{\longrightarrow}
\newcommand{\respect}{\rightsquigarrow}
\newcommand{\compatible}{\leftrightsquigarrow}
\newcommand{\upc}[1]{\overset {\lower 0.3ex \hbox{${\;}_{\circ}$}}{#1}}
\newcommand{\Gmlog}{\bG_{m, \log}}
\newcommand{\Gm}{\bG_m}
\newcommand{\ep}{\varepsilon}
\newcommand{\Spec}{\operatorname{Spec}}
\newcommand{\val}{{\mathrm{val}}} 
\newcommand{\n}{\operatorname{naive}}
\newcommand{\bs}{\operatorname{\backslash}}
\newcommand{\Gal}{\operatorname{{Gal}}}
\newcommand{\gal}{{\rm {Gal}}({\bar \Q}/{\Q})}
\newcommand{\galp}{{\rm {Gal}}({\bar \Q}_p/{\Q}_p)}
\newcommand{\gall}{{\rm{Gal}}({\bar \Q}_\ell/\Q_\ell)}
\newcommand{\wep}{W({\bar \Q}_p/\Q_p)}
\newcommand{\wel}{W({\bar \Q}_\ell/\Q_\ell)}
\newcommand{\Ad}{{\rm{Ad}}}
\newcommand{\BS}{{\rm {BS}}}
\newcommand{\even}{\operatorname{even}}
\newcommand{\End}{{\rm {End}}}
\newcommand{\odd}{\operatorname{odd}}
\newcommand{\GL}{\operatorname{GL}}
\newcommand{\np}{\text{non-$p$}}
\newcommand{\g}{{\gamma}}
\newcommand{\G}{{\Gamma}}
\newcommand{\Lam}{{\Lambda}}
\newcommand{\La}{{\Lambda}}
\newcommand{\lam}{{\lambda}}
\newcommand{\la}{{\lambda}}
\newcommand{\uL}{{{\hat {L}}^{\rm {ur}}}}
\newcommand{\uQp}{{{\hat \Q}_p}^{\text{ur}}}
\newcommand{\sel}{\operatorname{Sel}}
\newcommand{\dt}{{\rm{Det}}}
\newcommand{\Sig}{\Sigma}
\newcommand{\fil}{{\rm{fil}}}
\newcommand{\SL}{{\rm{SL}}}
\newcommand{\spl}{{\rm{spl}}}
\newcommand{\st}{{\rm{st}}}
\newcommand{\Isom}{{\rm {Isom}}}
\newcommand{\Mor}{{\rm {Mor}}}
\newcommand{\bg}{\bar{g}}
\newcommand{\id}{{\rm {id}}}
\newcommand{\cone}{{\rm {cone}}}
\newcommand{\al}{a}
\newcommand{\ChL}{{\cal{C}}(\La)}
\newcommand{\Image}{{\rm {Image}}}
\newcommand{\toric}{{\operatorname{toric}}}
\newcommand{\torus}{{\operatorname{torus}}}
\newcommand{\Aut}{{\rm {Aut}}}
\newcommand{\Qp}{{\mathbb{Q}}_p}
\newcommand{\barQp}{{\mathbb{Q}}_p}
\newcommand{\Qpur}{{\mathbb{Q}}_p^{\rm {ur}}}
\newcommand{\Zp}{{\mathbb{Z}}_p}
\newcommand{\Zl}{{\mathbb{Z}}_l}
\newcommand{\Ql}{{\mathbb{Q}}_l}
\newcommand{\Qlur}{{\mathbb{Q}}_l^{\rm {ur}}}
\newcommand{\F}{{\mathbb{F}}}
\newcommand{\eps}{{\epsilon}}
\newcommand{\epsLa}{{\epsilon}_{\La}}
\newcommand{\epsLaVxi}{{\epsilon}_{\La}(V, \xi)}
\newcommand{\epsOLaVxi}{{\epsilon}_{0,\La}(V, \xi)}
\newcommand{\Qplin}{{\mathbb{Q}}_p(\mu_{l^{\infty}})}
\newcommand{\otimesQplin}{\otimes_{\Qp}{\mathbb{Q}}_p(\mu_{l^{\infty}})}
\newcommand{\galFl}{{\rm{Gal}}({\bar {\Bbb F}}_\ell/{\Bbb F}_\ell)}
\newcommand{\gallur}{{\rm{Gal}}({\bar \Q}_\ell/\Q_\ell^{\rm {ur}})}
\newcommand{\galFF}{{\rm {Gal}}(F_{\infty}/F)}
\newcommand{\galFv}{{\rm {Gal}}(\bar{F}_v/F_v)}
\newcommand{\galF}{{\rm {Gal}}(\bar{F}/F)}
\newcommand{\epsVxi}{{\epsilon}(V, \xi)}
\newcommand{\epsOVxi}{{\epsilon}_0(V, \xi)}
\newcommand{\plim}{\lim_
{\scriptstyle 
\longleftarrow \atop \scriptstyle n}}
\newcommand{\sig}{{\sigma}}
\newcommand{\ga}{{\gamma}}
\newcommand{\del}{{\delta}}
\newcommand{\Vss}{V^{\rm {ss}}}
\newcommand{\Bst}{B_{\rm {st}}}
\newcommand{\Dpst}{D_{\rm {pst}}}
\newcommand{\Dcrys}{D_{\rm {crys}}}
\newcommand{\DdR}{D_{\rm {dR}}}
\newcommand{\Fin}{F_{\infty}}
\newcommand{\Kla}{K_{\lambda}}
\newcommand{\Ola}{O_{\lambda}}
\newcommand{\Mla}{M_{\lambda}}
\newcommand{\Det}{{\rm{Det}}}
\newcommand{\Sym}{{\rm{Sym}}}
\newcommand{\LaSa}{{\La_{S^*}}}
\newcommand{\cX}{{\cal {X}}}
\newcommand{\MHG}{{\frak {M}}_H(G)}
\newcommand{\tauMla}{\tau(M_{\lambda})}
\newcommand{\Fvur}{{F_v^{\rm {ur}}}}
\newcommand{\Lie}{{\rm {Lie}}}
\newcommand{\cL}{{\cal {L}}}
\newcommand{\cW}{{\cal {W}}}
\newcommand{\fq}{{\frak {q}}}
\newcommand{\cont}{{\rm {cont}}}
\newcommand{\SC}{{SC}}
\newcommand{\Om}{{\Omega}}
\newcommand{\dR}{{\rm {dR}}}
\newcommand{\crys}{{\rm {crys}}}
\newcommand{\hatSig}{{\hat{\Sigma}}}
\newcommand{\rdet}{{{\rm {det}}}}
\newcommand{\ord}{{{\rm {ord}}}}
\newcommand{\BdR}{{B_{\rm {dR}}}}
\newcommand{\BdRO}{{B^0_{\rm {dR}}}}
\newcommand{\Bcrys}{{B_{\rm {crys}}}}
\newcommand{\Qw}{{\mathbb{Q}}_w}
\newcommand{\barkappa}{{\bar{\kappa}}}
\newcommand{\cP}{{\Cal {P}}}
\newcommand{\cZ}{{\Cal {Z}}}
\newcommand{\oppLa}{{\Lambda^{\circ}}}

\begin{abstract}
 We extend SL(2)-orbit theorems for degeneration of mixed Hodge structures to a situation in which we do not assume the polarizability of graded quotients. We also obtain analogous results on Deligne systems.  
\end{abstract}
\renewcommand{\thefootnote}{\fnsymbol{footnote}}
\footnote[0]{Primary 14D07; Secondary 14C30, 58A14. 
The author is partially supported by an NSF grant.}

\section{Introduction} 

\begin{para}
In this paper, we show that the SL(2)-orbit theorems on the degeneration of Hodge structures (\cite{Scw},  \cite{CKS} and \cite{KNU1}) hold in  a situation in which we do not assume the polarizability of the graded quotients for the weight filtration. We also obtain analogous results on Deligne systems.

\end{para}

\begin{para}
Recall that a Deligne system of $n$ variables is $(V, W, N_1, \dots, N_n, \alpha)$ where $V$ is a finite dimensional vector space over a field $E$ of characteristic $0$, $W$ is a finite increasing filtration on $V$ (called the weight filtration), $N_1,\dots, N_n: V\to V$ are mutually commuting nilpotent linear operators (called the monodromy operators) which respect $W$, $\alpha$ is an action of the multiplicative group $\mathbb{G}_m$ on $V$, satisfying certain conditions (\cite{Scc};  see also \ref{defDH} of this paper for a review).

In this paper, we define a similar notion Deligne-Hodge system (DH system in short) of $n$ variables, which is $(V, W, N_1, \dots, N_n, F)$ where $(V, W, N_1,\dots, N_n)$ has the same properties as in the definition of Deligne system of $n$ variables with $E=\R$, and $F$ is a decreasing filtration on $V_\C=\C\otimes_{\R} V$ (called the Hodge filtration) satisfying certain conditions (see \ref{defDH}). 
A Deligne-Hodge system of zero variable is nothing but a mixed $\R$-Hodge structure. 

The notion Deligne-Hodge system is similar to the notion infinitesimal mixed Hodge module (IMHM) of Kashiwara (\cite{Kas}; see also \ref{IMHM} of this paper for a review). In fact, if $(V, W, N_1, \dots, N_n, F)$ is an IMHM, then it is a DH system of $n$ variables.
In the definition of DH system, we do not assume the polarizability of the graded quotients for weight filtration which was assumed for IMHM. Another difference is that in the definition of DH system, the order of $(N_1,\dots, N_n)$ matters though it does not matter for IMHM. 

\end{para}

\begin{para}

SL(2)-orbit theorems are statements on the
properties of $\exp(\sum_{j=1}^n iy_jN_j)F$ 
for an IMHM $(V, W, N_1,\dots, N_n, F)$ of $n$ variables in the situation  $y_j/y_{j+1}\to \infty$ ($1\leq j\leq n$, $y_{n+1}$ denotes $1$). 
  (\cite{Scw} treats the pure case with $n=1$, \cite{CKS} treats the pure case in general, and \cite{KNU1} generalizes it to the mixed case). 
In this paper, we prove the following Theorem \ref{thm1}  which 
shows that SL(2)-orbit theorems in \cite{Scw}, \cite{CKS}, \cite{KNU1} 
are generalized to DH systems.   

\end{para}

\begin{thm}\label{thm1} Let $(V, W, N_1, \dots, N_n, F)$ be a DH system of $n$ variables. Then for  $N'_j=\sum_{k=1}^j a_{j,k}   N_k$ ($1\leq j\leq n$) with $a_{j,k}>0$ ($1\leq k \leq  j\leq n$) such that $a_{j,k}/a_{j,k+1}\gg 0$ ($1\leq k <j\leq n$), 
 $(V, W, N'_1, \dots, N'_n, F)$
is an IMHM of $n$ variables.

\end{thm}

For example, if $(V, W, N_1, N_2, F)$ is a DH system of two variables, $(V,W, N_1, aN_1+N_2, F)$ for $a\gg 0$ is an IMHM. 

For  $N'_j$ as in Theorem \ref{thm1}, if  $y_j/y_{j+1}\to \infty$ ($1\leq j\leq n$),  we have
$\sum_{j=1}^n y_jN_j=\sum_{j=1}^\infty y'_jN'_j$ with $y'_j/y'_{j+1}\to \infty$ ($1\leq j\leq n$). 
 Hence the property of $\exp(\sum_{j=1}^n iy_jN_j)F$ in the situation $y_j/y_{j+1}\to \infty$ ($1\leq j\leq n$)
 for a DH system
 is reduced to the case of IMHM.

\begin{para}

We have a canonical functor from the category of DH systems of $n$ variables to the category of Deligne systems of $n$ variables over $\R$, which has the shape $(V, W, N_1, \dots, N_n, F)\mapsto (V, W, N_1,\dots, N_n, \alpha)$ for a canonically defined $\alpha$ (see Section 2.2). We have also a canonical functor from the category of Deligne systems of $n$ variables over $\R$ or over $\C$ to the category of DH systems of $n$ variables which has the shape
$(V, W, N_1, \dots, N_n, \alpha)\mapsto (V^{\oplus 2}, W^{\oplus 2}, N_1^{\oplus 2}, \dots, N_n^{\oplus 2}, F)$ for a canonically defined $F$ (Section 2.3). Here in the case of Deligne system over $\C$, $V^{\oplus 2}$ is regarded as an $\R$-vector space by restriction of scalers. We study Deligne systems and DH systems by using these two functors applying the results on one to the other.

From the above theorem on DH systems, we obtain the following theorem on Deligne systems. 
\end{para}

\begin{thm}\label{thm2}

 Let $(V, W, N_1, \dots, N_n, \alpha)$ be a Deligne system
  of $n$ variables  over $\R$ or over $\C$. Then for  $N'_j=\sum_{k=1}^j a_{j,k}   N_k$ with $a_{j,k}>0$ ($1\leq k \leq  j\leq n$) such that $a_{j,k}/a_{j,k+1}\gg 0$ ($1\leq k <j\leq n$), 
   the DH system $(V^{\oplus 2}, W^{\oplus 2}, (N'_1)^{\oplus 2}, \dots, (N'_n)^{\oplus 2}, F)$ of $n$ variables
   associated to the Deligne system
   $(V^{\oplus 2}, W^{\oplus 2}, (N'_1)^{\oplus 2}, \dots, (N'_n)^{\oplus 2}, \alpha^{\oplus 2})$
    is an IMHM. 

\end{thm}

This tells that, roughly speaking, any Deligne system of $n$ variables underlies some IMHM if it is modified in an elementary way  and replaced by the direct sum of two copies of it.

From Theorem \ref{thm1} (resp. Theorem \ref{thm2}) and the SL(2)-orbit theorem in \cite{KNU1}, we have the part on DH systems  (resp. Deligne systems) in the following theorem. 

\begin{thm}\label{Dthm} (1) Let $(V, W, N_1, \dots, N_n, F)$ be a DH system of $n$ variables. 
If  $y_j/y_{j+1}\gg 0$ ($1\leq j\leq n$, $y_{n+1}$ denotes $1$), $(V, W, \exp(\sum_{j=1}^n iy_jN_j)F)$ is a mixed Hodge structure. 
 The splitting of $W$ associated to this mixed Hodge structure (\ref{delta}) converges when  
 $y_j/y_{j+1}\to \infty$ ($1\leq j\leq n$).

\medskip

(2)
Let  $E$ be $\R$ or $\C$, and let $(V, W, N_1, \dots, N_n, \alpha)$ be a Deligne system of $n$ variables over $E$. 
Let  $W'$ be  the increasing filtration on $V$ defined by $\alpha$ (for $w\in \Z$,
$W'_w$ is defined as the sum of the weight $k$ part for $\alpha$ for all $k\leq w$). Then if  $y_j> 0$ ($1\leq j\leq n$) and  $y_j/y_{j+1}\gg 0$ ($1\leq j<n$), $W'$ is the relative monodromy filtration (\ref{relmo}) of $\sum_{j=1}^n y_jN_j$ with respect to $W$. The splitting of $W$ defined by the  
Deligne system $(V, W, \sum_{j=1}^n y_jN_j, \alpha)$ of one variable (\ref{Dthm1}) converges when $y_j>0$ ($1\leq j\leq n$) and $y_j/y_{j+1}\to \infty$ ($1\leq j<n$).

\end{thm}

In Theorem \ref{Wspas} in Section 4.2,  we will give more precise descriptions of the convergences in (1) and (2) of this theorem.

The following result is  deduced from Theorem \ref{thm1} and Theorem \ref{thm2} and from the fact that the category of IMHM of $n$ variables is an abelian category (\cite{Kas}). 

\begin{thm}\label{thm3} The category of Deligne systems of $n$ variables over a field $E$ of characteristic $0$ is an abelian category. The category of DH systems of $n$ variables is an abelian category. In these categories, the underlying vector space $V$ of the kernel (resp. cokernel) of a morphism $A\to B$ is the kernel (resp. cokernel) of the map of the underlying spaces,  and $W$, $N_j$, etc. of the kernel (resp. cokernel) are the ones induced from those of $A$ (resp. $B$). 

\end{thm}

\begin{para}
We expect that  results of this paper are useful to generalize the works \cite{KNU} on classifying spaces of degenerating Hodge structures to a situation where we do not assume the polarizability of the graded quotients for the weight filtration. We also expect that the study on Deligne systems as in this paper are useful in the studies (like \cite{BK} and 
\cite{KK})  which treat degeneration of motives over non-archimedean local fields.

The author thanks Spencer Bloch, Chikara Nakayama, and Sampei Usui for the joint works with the author (\cite{BK}, \cite{KNU1}, \cite{KNU}) which inspired this work.  
\end{para}

\section{Deligne systems and Deligne-Hodge systems}

\subsection{Definitions}

\begin{sbpara}\label{relmo}
We first review the notion relative monodromy filtration defined by Deligne  (\cite{De} 1.6.13).

Let $V$ be an abelian group, let $W=(W_w)_{w\in \bZ}$ be a finite increasing  filtration on $V$, and let $N:V\to V$ be a nilpotent homomorphism such that $NW_w\subset W_w$ for alll $w\in \Z$.  

Then a finite  increasing filtration $W'=(W'_w)_{w\in \bZ}$ on $V$ is called the 
{\it relative monodromy filtration of $N$ with respect to $W$} if it satisfies the following conditions (i) and (ii). 

(i) $NW'_w\subset W'_{w-2}$ for any $w\in \bZ$. 

 (ii) For any $w\in \bZ$ and $m\geq 0$, we have an isomorphism $N^m: \gr^{W'}_{w+m}\gr^W_w\overset{\cong}\longrightarrow  \gr^{W'}_{w-m}\gr^W_w.$

The relative monodromy filtration of $N$ with respect to $W$ need not exist. If it exists, it is unique (\cite{De} 1.6.13).

If $V$ is a vector space over a field $E$ and if $W_w$ are $E$-linear subspaces and $N$ is $E$-linear, the realtive monodromy filtration consists of $E$-linear subspaces of $V$ if it exists. 

\end{sbpara}

\begin{sbpara}\label{defDH} We review the notion Deligne system of $n$ variables (\cite{Scc}), and define the notion Deligne-Hodge system of $n$ variables.

A Deligne system 
 over a field $E$ of characteristic $0$ (resp. A Deligne-Hodge system)  is $$(V, W,  N_1, \dots, N_n, \alpha)\quad (\text{resp.}\;\;
 (V, W, N_1,\dots, N_n, F))$$
where 

$V$ is a finite dimensional $E$ (resp. $\R$)-vector space,

$W$ is a finite increasing  filtration on $V$,

$N_j$ are  linear operators $V\to V$,

$\alpha$ is an action of $\mathbb{G}_m$ on $V$ (resp. $F$ is a finite decreasing filtration on $V_{\C}=\C\otimes_{\R} V$),

\noindent
satisfying the following conditions (a), (b), (c), (d), and (e) (resp. (a), (b), (c), (d) and (f1), (f2)).

(a) The operators $N_1,\dots, N_n$ are nilpotent,  mutually commute, and respect $W$.

(b) There are finite increasing filtrations $W^{(j)}$ ($0\leq j\leq n$) such that $W^{(0)}=W$ and such that for $1\leq j\leq n$, $W^{(j)}$ is the relative monodromy filtration of $N_j$ with respect to $W^{(j-1)}$. 

(c) Let $1\leq j\leq n$, $0\leq k<j-1$, $w\in \Z$,   and let  $U=W^{(k)}_w$. Then  the restriction $W^{(j)}|_U$ of $W^{(j)}$ to $U$ is the relative monodromy filtration of $N_j|_U$ with respect to $W^{(j-1)}|_U$. 

(d) $N_j(W^{(k)}_w)\subset W^{(k)}_w$ for any $j, k, w$, and $N_j(W^{(k)}_w)\subset W^{(k)}_{w-2}$ if $k\geq j$,

(e) $\alpha$ splits $W^{(n)}$,  $W^{(j)}_w$ is stable under the action $\alpha$ of $\mathbb{G}_m$ for any $0\leq j<n$ and $w\in \Z$, and $N_j$ is of weight $-2$ for $\alpha$ (that is, $\alpha(a)N_j\alpha(a)^{-1}=a^{-2}N_j$ for any $a$) for any $1\leq  j\leq n$. 

(f1) $N_jF^p\subset F^{p-1}$ for any $1\leq j\leq n$ and $p\in \Z$.

(f2) $(W^{(n)}, F)$ is a mixed Hodge structure. Furthermore, for $1\leq k< n$, $w\in \Z$ and for $U=W^{(k)}_w$, $(W^{(n)}|_U, F|_U)$ is a mixed Hodge structure. 

\end{sbpara}

\begin{sbpara}
We denote the category of Deligne systems of $n$ variables over $E$ by $D_{n,E}$.

We denote the category of Deligne-Hodge systems of $n$ variables by $DH_n$.

\end{sbpara}

\begin{sbpara}
For example, a Deligne system of zero variable over $E$ is nothing but a finite dimensional $E$-vector space endowed with an action of $\mathbb{G}_m$. 

A Deligne-Hodge system of zero variable is just a mixed $\R$-Hodge structure. In this paper, we call a mixed $\R$-Hodge structure just a mixed Hodge structure.

\end{sbpara}

\begin{sbpara}
A Deligne system of  one variable over $E$ is nothing but $(V, W,N, \alpha)$ where 

$V$ is a finite dimensional $E$-vector space, 

$W$ is a finite increasing filtration on $V$, 

$N$ is a nilpotent linear map $V\to V$ such that $N(W_w)\subset W_w$ for any $w\in \Z$, 

$\alpha$ is an action of $\mathbb{G}_m$ on $V$,

\noindent
such that  $W_w$ is stable under the action $\alpha$ of $\mathbb{G}_m$ for any $w\in \Z$, $N$ is of weight $-2$ for $\alpha$, and if we define $W'_w$ to be the sum of the weight $k$ part of $\alpha$ for all $k\leq w$, then $W'$ is the relative monodromy filtration of $N$ with respect to $W$.

\end{sbpara} 

\begin{sbpara} Both the categories $D_{n,E}$ and $DH_n$ have direct sum, tensor products, symmetric powers, exterior powers, duals, and Tate twists, defined in the evident manners.

\end{sbpara}

The following is easy to see.
\begin{sblem}\label{EandE'} Let $E$ be a field of characteristic $0$ and let $H=(V, W, N_1, \dots, N_n, \alpha)$ be as in the hypothesis of the definition of the notion Deligne system of $n$ variables (we do not assume (a)--(e)). 

(1) Let $E'$ be a field which contains $E$ as a subfield and let $H'=E'\otimes_E H$. Then $H$ is in $D_{n,E}$ if and only if $H'$ is in $D_{n,E'}$. 

(2) Let $E'$ be a subfield of $E$ such that $E$ is a finite extension of $E'$. Let $H'$ be $H$ but $V$ in $H'$ is regarded as an $E'$-vector space by the restriction of scalers. Then $H$ is in $D_{n,E}$ if and only of $H'$ is  in $D_{n,E'}$.
\end{sblem}

The following is also easy to see.

\begin{sblem}  Let $(V, W, N_1,\dots, N_n, \alpha)$ (resp. $(V, W, N_1,\dots, N_n, F)$)
be an object of $D_{n,E}$ (resp. $DH_n$). Then for any  $a_{j,k}\in E$ (resp. $\R$) ($1\leq k\leq j\leq n$) such that $a_{j,j}\neq 0$ ($1\leq j\leq n$), 
if we put $N'_j=\sum_{k=1}^j a_{j, k}N_k$ for $1\leq j\leq n$, then $(V, W, N'_1, \dots, N'_n, \alpha)$ (resp.
 $(V, W, N'_1, \dots, N'_n, F)$
belongs to $D_{n,E}$ (resp. $DH_n$).

\end{sblem}

\begin{sbpara}\label{IMHM}The notion DH system of $n$ variables is similar to the notion  IMHM of Kashiwara. We review the notion IMHM (in fact we consider in this paper only IMHM which has $\R$-structure, and we call such IMHM just IMHM in this paper).

An IMHM of $n$ variables is $(V, W, N_1, \dots, N_n, F)$ as in the hypothesis of the definition of the notion DH system of $n$ variables, satisfying the following conditions (a), (f1), (g), (h). 

(a) The same as (a) in \ref{defDH}.

(f1) The same as (f1) in \ref{defDH}. 

(g) For each $w\in \Z$, there is a non-degenerate $\R$-bilinear form $\langle\;,\;\rangle : \gr^W_w\times \gr^W_w\to \R$ which is symmetric if $w$ is even and anti-symmetric if $w$ is odd such that 
$\langle N_ju, v\rangle_w+\langle u, N_jv\rangle_w=0$ for any $j$ and any $u, v\in \gr^W_w$ and such that if  $y_j\gg 0$ ($1\leq j\leq n$), $(\gr^W_w, (\;,\;), \exp(\sum_{j=1}^n iy_jN_j)F(\gr^W_w))$ is a polarized Hodge structure. Here $F(\gr^W_w)$ denotes the filtration on $\gr^W_{w,\C}$ induced by $F$.

(h) For $1\leq j\leq n$, the relative monodromy filtration of $N_j$ with respect to $W$ exists. 
\end{sbpara}

By the arguments in \cite{Scc}, Section 3, Example 2, we have:

\begin{sbprop} An
IMHM of $n$ variables is a DH system of $n$ variables. 

\end{sbprop}

\subsection{A functor $DH_n\to D_{n,\R}$}

We define a functor $DH_n\to D_{n,\R}$.

\begin{sbpara}\label{delta}

 We review that for a mixed Hodge structure $(V, W, F)$, we have a canonical splitting of $W$.  (This canonical splitting is called the SL(2)-splitting in \cite{BP}).

 There is a unique pair $(s', \delta)$ of a splitting $s': \gr^W=\oplus_{w\in \Z} \; \gr^W_w\overset{\cong}\to V$ of $W$ and a linear map $\delta:\gr^W\to \gr^W$  such that the  Hodge $(p,q)$-component $\delta_{p,q}$ of $\delta$ for $F(\gr^W)$ is zero unless $p<0$ and $q<0$ and such that  $F=s'(\exp(i\delta)\gr^W(F))$. (See  \cite{CKS} 2.20.)

The canonical splitting $s$ of $W$ is a modification of this $s'$. It is defined by $s=s'\circ \exp(-\zeta)$ where 
$\zeta: \gr^W\to \gr^W$ is the linear map  which is determined by $\delta$ as Lie polynomial of $\delta_{p,q}$ as in \cite{CKS} 6.60. 

Any morphisms of mixed Hodge structures commutes with canonical splittings. 
\end{sbpara}

\begin{sblem}
Let  $(V, W, N_1, \dots, N_n, F)$ be a DH system of $n$ variables, and let $\alpha$ be the canonical splitting of $W^{(n)}$ associated to the mixed Hodge structure $(W^{(n)}, F)$.
Then $(V, W, N_1, \dots, N_n, \alpha)$ is a Deligne system of $n$ variables.
\end{sblem}

\begin{pf} It is sufficient to prove the following (1) and (2).

(1) For any $0\leq j<n$ and $w\in \Z$, $W^{(j)}_w$ is stable under the action $\alpha$ of $\mathbb{G}_m$.

(2) For any $1\leq j\leq n$, $N_j$ is of weight $-2$ for $\alpha$.

We prove (1). Let $U=W_w^{(j)}$. The inclusion map $U\to V$ is a morphism of mixed Hodge structures $(W^{(n)}|_U, F|_U)\to (W^{(n)}, F)$. Hence the canonical splitting of 
$W^{(n)}|_U$ associated to the mixed Hodge structure $(W^{(n)}|_U, F|_U)$ and the canonical splitting of $W^{(n)}$ associated to the mixed Hodge structure 
$(W^{(n)}, F)$ (that is $\alpha$) is compatible. This proves (1). 

We prove (2). By $N_jF^p\subset F^{p-1}$ for any $p$, $N_j$ is a morphism of mixed Hodge structures $(W^{(n)}, F)\to (W^{(n)}(-1), F(-1))$ where $(-1)$ is the Tate module. Hence via $N_j$,  the canonical splitting of $W^{(n)}$ associated to the mixed Hodge structure $(W^{(n)}, F)$ is compatible with the canonical splitting of $W^{(n)}(-1)$ associated to the mixed Hodge structure $(W^{(n)}(-1), F(-1))$. This proves (2).
\end{pf}

\begin{sbpara}  Thus we obtained  the functor
$$DH_n\to D_{n,\R}\;;
 (V, W, N_1, \dots, N_n, F)\mapsto (V, W, N_1, \dots, N_n, \alpha).$$

\end{sbpara}

\subsection{A functor $D_{n,E}\to DH_n$ for $E=\R$ or $\C$}

\begin{sbpara}
For  $E=\R$ or $\C$, we define a functor $D_{n,E}\to DH_n$. This functor in the case $E=\C$ is defined to be the composition $$D_{n,\C}\to D_{n,\R}\to DH_n$$ where the first functor is to regard a $\C$-vector space as an $\R$-vector space.  So, we assume $E=\R$. 

\end{sbpara}

\begin{sbpara}

Let $(V, W, N_1, \dots, N_n, \alpha)$ be a Deligne system of $n$ variables over $\R$. We define a decreasing filtration $F$ on $V^{\oplus 2}_\C$ as follows. For $w\in \Z$, let $V_w$ be the weight $w$ part of $V$ with respect to the action $\alpha$ of $\mathbb{G}_m$. We define $F$ as a direct sum of decreasing filtrations on $V_{w,\C}^{\oplus 2}$. If $w$ is an even integer $2r$, define the filtration $F$ on $V_{w,\C}^{\oplus 2}$ by $F^r=V_{w,\C}^{\oplus 2}$ and $F^{r+1}=0$. 
 If $w$ is an odd integer $2r+1$, define the filtration $F$ on $V_{w,\C}^{\oplus 2}$ as follows: $F^r=V_{w,\C}^{\oplus 2}$, $F^{r+2}=0$, and $F^{r+1}$ is 
the $\C$-subspace of $V_{w,\C}^{\oplus 2}$ generated by $(i\otimes x, 1\otimes x)$ $(x\in V_w)$.

\end{sbpara}

\begin{sblem} This $(V^{\oplus 2}, W^{\oplus 2}, N_1^{\oplus 2}, \dots, N_n^{\oplus 2}, F)$  is a DH system of $n$ variables. 

\end{sblem}

This is checked easily.

\begin{sbpara}  Thus we obtained  the functor
$$D_{n,\R}\to DH_n\;\;;\;\;
 (V, W, N_1, \dots, N_n, \alpha)\mapsto (V^{\oplus 2}, W^{\oplus 2}, N_1^{\oplus 2}, \dots, N_n^{\oplus 2}, F).$$
 The composition $D_{n,\R}\to DH_n\to D_{n,\R}$ with the functor in Section 2.2 is
 $$ (V, W, N_1, \dots, N_n, \alpha)\mapsto (V^{\oplus 2}, W^{\oplus 2}, N_1^{\oplus 2}, \dots, N_n^{\oplus 2}, \alpha^{\oplus 2}).$$
 \end{sbpara}

\section{SL(2)-orbits}

\subsection{Splittings of Deligne}

We review two theorems of Deligne on splittings of weight filtrations of Deligne systems which are introduced in \cite{Scc} as Theorem 1 and Theorem 2, respectively. 

\begin{sbpara}\label{3.1.1}  First we review the notion primitive component. Let $V$ be an abelian group, let $W$ be a finite increasing filtration on $V$, and let $N: V\to V$ be a nilpotent endomorphism which respects $W$. Assume that the relative monodromy filtration $W'$ of $N$ with respect to $W$ exists. 
Let $w\in \Z$, $m\geq 0$. Then $\gr^{W'}_{w+m}\gr^W_w=A \oplus B$, where $A$ is the kernel of $\gr^{W'}_{w+m}\gr^W_w\overset{N^{m+1}}\longrightarrow \gr^{W'}_{w-m-2}\gr^W_w$ and $B$ is the image of $N:\gr^{W'}_{w+m+2}\gr^W_w\to \gr^{W'}_{w+m}\gr^W_w$. The component $A$ is called the {\it primitive component} of $\gr^{W'}_{w+m}\gr^W_w$.

\end{sbpara}

\begin{sbpara}\label{3.1.2} Let $V, W, N, W'$ be as in \ref{3.1.1}. Denote the filtration on $Hom(V, V)$ induced by $W$ (resp. $W'$) by $W_{\bullet}Hom(V, V)$ (resp. $W'_{\bullet}Hom(V, V)$). Then $W'_{\bullet}Hom(V, V)$ is the relative monodromy filtration of the nilpotent homomorphism $Ad(N): Hom(V, V)\to Hom(V, V)$ with respect to $W_{\bullet}Hom(V, V)$. 

\end{sbpara}

\begin{sbpara} \label{Dthm1} Let $(V, W, N, \alpha)$ be a Deligne system of one variable. The first theorem of Deligne is that there is a unique action $\tau=(\tau_0,\tau_1)$ of $\mathbb{G}_m^2$ on $V$ satisfying the following conditions (i)--(iii).

(i) $\tau_1=\alpha$.

(ii) $\tau_0$ splits $W^{(0)}=W$.

(iii) For $k\geq 1$, let $N_{-k}\in \gr^W_{-k}\Hom(V, V)$ be the weight $-k$ part of $N$ with respect to the action $\tau_0$ of $\mathbb{G}_m$ on $V$. Then $N_{-1}=0$, and for any $k\geq 2$, the class of $N_{-k}$ in $\gr^{W'}_{-2}\gr^W_{-k}\Hom(V, V)$ belongs to the primitive component.

\end{sbpara}

\begin{sbpara}\label{Dthm2}
The second theorem of Deligne is the following.

Let $(V, W, N_1,\dots, N_n, \alpha)$ be a Deligne system of $n$ variables. Then there is a unique action of $\tau=(\tau_j)_{0\leq j\leq n}$ of $\mathbb{G}_m^{n+1}$ on $V$ satisfying the following conditions (i) and (ii).

\medskip
(i) $\tau_n=\alpha$.

(ii) For $1\leq j\leq n$, $(V, W^{(j-1)}, N_j, \tau_j)$ is a Deligne system of one variable, and the action $(\tau_{j-1}, \tau_j)$ of $\mathbb{G}_m^2$ coincides with the action of $\mathbb{G}_m^2$ in \ref{Dthm1} associated to this Deligne system of one variable. 

\medskip

Furthermore, for this $\tau$, we have the following (iii), (iv), (v).

(iii) For $0\leq j\leq n$, $\tau_j$ splits $W^{(j)}$.

(iv) For $1\leq j\leq k\leq n$, $N_j$ is of weight $-2$ for $\tau_k$. 

(v) Let $1\leq j\leq n$, and let $\hat N_j$ be the component of $N_j$ of weight $0$ for $\tau_{j-1}$. Then $\hat N_j$ is of weight $0$ for $\tau_k$ for any $0\leq k<j$.

\end{sbpara}

\begin{sbpara}
If $(V, W, N_1, \dots, N_n, F)$ is a DH system of $n$ variables, we have the associated action $\tau=(\tau_j)_{0\leq j\leq n}$ of $\mathbb{G}_m^{n+1}$ on $V$ defined by the corresponding Deligne system $(V,W, N_1, \dots, N_n, \alpha)$ (Section 2.2). 
\end{sbpara}

\subsection{SL(2)-orbits}

\begin{sbpara}

We say a Deligne system $(V, W, N_1,\dots, N_n,\alpha)$ of $n$ variables is  an SL(2)-orbit if 
$$\tau_k(a)N_j\tau_k(a)^{-1}=N_j\quad\text{ for}\;  0\leq k<j\leq n$$
for any $a$ (that is, $N_j$ is of weight $0$ for $\tau_k$ for $0\leq k<j\leq n$). We denote the full subcategory of $D_{n,E}$ consisting of SL(2)-orbits by $\hat D_{n,E}$.

We say a DH system $(V, N_1, \dots, N_n, F)$ of $n$ variables is an SL(2)-orbit if 
$$\tau_k(a)N_j\tau_k(a)^{-1}=N_j\quad\text{ for  $0\leq k<j\leq n$\;\; and}\quad \tau_k(a)F=F\quad \text{for $0\leq k \leq n$}$$
for any $a$. 
 We denote  the full subcategory of $DH_n$ consisting of SL(2)-orbits by $\hat DH_n$. 
\end{sbpara}

 \begin{sblem} In the category $D_{n,E}$ (resp. $DH_n$), $\hat D_{n,E}$ (resp. $\hat DH_n$) is  stable under taking direct sums, tensor products, symmetric powers, exterior powers, duals, and Tate twists. 
  \end{sblem}

\begin{sbpara}\label{sleq1}
As is easily seen, we have the following  equivalence of  categories between $\hat D_{n,E}$ and the category of finite dimensional representations of $\mathbb{G}_m\times SL(2)^n$ over $E$:  For an object $(V, W, N_1, \dots, N_n, \alpha)$  of $D_{n,E}$ with the associated $(\tau_j)_{0\leq j\leq n}$,  the  corresponding representation is $(V, \rho)$ where $\rho$ is the action of
$\mathbb{G}_m\times SL(2)^n$ on $V$ characterized by the following properties (i), (ii), (iii).

(i) The action of $\mathbb{G}_m =\mathbb{G}_m\times \{1\}\subset \mathbb{G}_m\times SL(2)^n$ is $\tau_0$.

(ii) For $1\leq j\leq n$, the action of $\begin{pmatrix} 1/a&0\\0&a\end{pmatrix}$ in the $j$-th SL(2) is $\tau_j(a)/\tau_{j-1}(a)$.

(iii)  In the action of ${\frak s}l(2)$ induced by the action of the $j$-th SL(2),
 $\begin{pmatrix} 0&1\\0& 0\end{pmatrix}\in {\frak s}l(2)$ acts as $N_j$. 

We have 

(iv) For $0\leq j\leq n$, $\tau_j(a)=\rho(a, g)$ where $g= (g_k)_{1\leq k\leq n}\in  SL(2)^n$ with $g_k=\begin{pmatrix} 1/a&1\\0&a\end{pmatrix}$ if $k\leq j$, and $g_k=1$ if $k>j$, 

Conversely, for a finite dimensional representation $(V, \rho)$ of $\mathbb{G}_m\times SL(2)^n$, the corresponding object $(V,W, N_1,\dots, N_n, \alpha)$ is given as follows. $W$ is defined by $\tau_0$, $N_j$ are given by the above (iii), $\alpha=\tau_n$ is given by the case $j=n$ of the above (iv). 
\end{sbpara}

\begin{sbpara}\label{[rho]} We consider $\hat DH_n$. For a finite dimensional representation $(V, \rho)$ of $\mathbb{G}_m\times SL(2)^n$ over $\R$ such that the action $\tau_n$ of $\mathbb{G}_m$ on $V$ defined by the case $j=n$ of \ref{sleq1} (iv) has only even weights, we have an object  $[\rho]$ of $\hat DH_n$ defined as follows. Let $(V, N_1, \dots, N_n, \alpha)$ be the object of $\hat D_{n,\R}$ corresponding to $(V, \rho)$ (so $\alpha=\tau_n$ has only even weights), and let $V_{2r}$ ($r\in \Z$) be the weight $2r$ part of $V$ with respect to $\alpha$.  Let $[\rho]=(V, W, N_1,\dots, N_n, F)$ where $F$ is the direct sum of the decreasing filtrations on $V_{2r,\C}$  defined by $F^rV_{2r,\C}=V_{2r,\C}$ and $F^{r+1}V_{2r,\C}=0$. 

Any objects of $\hat DH_n$ is isomorphic to a direct sum of objects of the form $[\rho] \otimes H$ where $H$ is a pure Hodge structure which we regard as an object of $\hat DH_n$ in the evident way: $N_j=0$ on $H$ for all $j$, and $W$ is  pure of weight the weight of $H$. More precisely, we have the following description \ref{eqhat} (2) of $\hat DH_n$. 

\end{sbpara}

The following \ref{eqhat} (1) is a consequence of \ref{sleq1} and the well known classification of representations of $SL(2)^n$. (2) is deduced  from (1) by using the functor $DH_n\to D_{n,\R}$ in Section 2.2. 

\begin{sbprop}\label{eqhat}
(1) For $1\leq j\leq n$, let $P_j$ be the object of $\hat D_{n,E}$ corresponding to the two dimensional representation of $\mathbb{G}_m\times SL(2)^n$ given by the projection to the $j$-th SL(2).  For $k\in \Z$, let $S_k$ be the object of $\hat D_{n,E}$ corresponding to the one dimensional representation of $\mathbb{G}_m\times SL(2)^n$ defined as $(a,g)\mapsto a^k$ ($a\in \mathbb{G}_m$, $g\in SL(2)^n$).

Then the category $\hat D_{n,E}$ is equivalent to the category of families $(H_{m,k})_{m\in \N^n,k\in \Z}$, where $H_{m,k}$ is a finite dimensional $E$-vector space for each $m, k$, satisfying $H_{m,k}=0$ for almost all $(m, k)$. The functor from the latter category to the former category  
$$(H_{m,k})_{m,k}\mapsto \oplus_{m,k}\; Sym^{m(1)}(P_1)\otimes \dots \otimes Sym^{m(n)}(P_n) \otimes S_k \otimes H_{m,k}.$$
gives an equivalence of categories. Here $H_{m,k}$ is regarded as an object of $\hat D_n$ in the following simple way. $V=H_{m,k}$, $W_0=V$, $W_{-1}=0$, $N_j=0$ for all $j$, $\alpha$ is trivial.

The inverse functor sends an object $(V, W, N_1, \dots, N_n, \alpha)$ to $(H_{m,k})_{m,k}$, where 
$$H_{m,k}=\{x\in V\; |\; N_j(x)=0 \; (1\leq j\leq n), \tau_j(a)x= a^k \prod_{\ell=1}^j a^{-m(\ell)}x\; (1\leq j\leq n)\;\text{for  any $a$}\}.$$

\medskip

(2) For $1\leq j\leq n$, let $P_j$ be the object $[\rho]$ of $\hat DH_n$ corresponding to the two dimensional representation $\rho$ of $\mathbb{G}_m\times SL(2)^n$ given by $(a,g)\mapsto ag_j$ ($a\in \mathbb{G}_m$, $g=(g_k)_k\in SL(2)^n$). 

Then the category $\hat DH_n$ is equivalent to the category of families $(H_{m,k})_{m\in \N^n,k\in \Z}$ where $H_{m,k}$ is a pure Hodge structure of weight $k$  for each $m, k$ satisfying $H_{m,k}=0$ for almost all $(m, k)$. The functor from the latter category to the former category  
$$(H_{m,k})_{m,k}\mapsto \oplus_{m,k}\; Sym^{m(1)}(P_1)\otimes \dots \otimes Sym^{m(n)}(P_n)  \otimes H_{m,k}.$$
gives an equivalence of categories. Here $H_{m,k}$ is regarded as an object of $\hat DH_n$ in the simple way explained in \ref{[rho]}.

The inverse functor sends an object $(V, W, N_1, \dots, N_n, F)$ to $(H_{m,k})_{m,k}$, where 
$$H_{m,k}=\{x\in V\; |\; N_j(x)=0 \; (1\leq j\leq n), \tau_j(a)\tau_{j-1}(a)^{-1}x= a^{-m(j)}x\; (1\leq j\leq n)\;\text{for any $a$}\}$$ 
whose Hodge filtration is  the restriction of $F$.

\end{sbprop}

\begin{sbprop}\label{sl2nm}

Let $(V, W, N_1, \dots, N_n, \alpha)$ be an object of  $\hat D_{n,E}$,  Fix $0\leq \ell\leq  j<k\leq n$. Then for any non-zero elements $y_t$ ($\ell+1\leq t\leq k$) of $E$, 
$W^{(k)}$ is the relative monodromy filtration of $\sum_{t=\ell+1}^k y_tN_t$ with respect to $W^{(j)}$. In other words, 
$(V, W^{(j)}, \sum_{t=\ell+1}^k   y_tN_t, \tau_k)$ is a Deligne system of one variable. 
\end{sbprop}

\begin{pf} By \ref{eqhat} (1), it is sufficient to check this in the cases of the objects $P_s$ ($1\leq s \leq n$) and $S_w$ ($w\in \Z$)  in \ref{eqhat} (1). These are checked easily. 
\end{pf}

\begin{sbprop}\label{sl2pol}  An object of $\hat DH_n$ is an IMHM. If $H=(V, W, N_1, \dots, N_n, F)$ is an object of  $\hat DH_n$, then for each $w\in \Z$, there is a non-degenerate  $\R$-bilinear form
$\langle\;,\rangle_w$ on $\gr^W_w$ such that for any $y_j>0$ ($1\leq j\leq n$), $(\gr^W_w, \langle \;,\;\rangle_w, \exp(\sum_{j=1}^n iy_jN_j)F(\gr^W_w))$ is a polarized Hodge structure
of weight $w$ and such that 
$$\langle \tau_j(a)u, \tau_j(a)v\rangle_w =a^{2w}\langle u, v\rangle_w,\quad \langle N_ju, v\rangle_w+\langle u, N_jv\rangle_w=0$$  for any $u,v\in \gr^W_w$  and for any $a$ and $j$.

\end{sbprop}

\begin{pf} By \ref{eqhat} (2), it is sufficient to prove this for the objects $P_j$ of $\hat DH_n$ ($1\leq j\leq n$) in \ref{eqhat} (2) and for a mixed Hodge structure regarded as an object of $\hat DH_n$ in the simple way in \ref{[rho]}. The case of mixed Hodge structure is well known (note that we are considering mixed $\R$-Hodge structures, so all mixed Hodge conjectures are polarizable as is well known). We consider the case of $P_j$. It is a two dimensional vector space $V$ over $\R$ with basis $(e_1, e_2)$, $W_1=V$, $W_0=0$, $N_je_2=e_1$, $N_je_1=0$, $N_k=0$ for any $k\neq j$, and $\alpha(a)e_1=e_1$, $\alpha(a)e_2=a^2e_2$. The condition (g) in \ref{IMHM} is satisfied because the anti-symmetric bilinear form  on $\gr^W_1$ defined by $\langle e_2, e_1\rangle_{1}=1$ satisfies the condition (g). The condition (h) in \ref{IMHM} is satisfied because $W^{(n)}$ is the relative monodromy filtration of $N_j$ with respect to $W$ and for $k\neq j$, $W$ is the relative monodromy filtration of $N_k=0$ with respect to $W$. 
\end{pf}

\subsection{Associated SL(2)-orbits}

\begin{sbpara}\label{hatF1}
For an object $(V, W, N_1,\dots, N_n, F)$ of $DH_n$, let  $$\hat F=s(F(\gr^{W^{(n)}}))$$ where $s:\gr^W\overset{\cong}\to V$ is the canonical splitting of $W^{(n)}$ associated to the mixed Hodge structure $(W^{(n)}, F)$.
\end{sbpara}

\begin{sbprop}\label{assl}
(1) Let $(V, W, N_1,\dots, N_n, \alpha)$ be an object of $D_{n,E}$. Then $(V, W, \hat N_1, \dots, \hat N_n,\alpha)$, where $\hat N_j$ are as in \ref{Dthm2}, is an object of $\hat D_{n,E}$. 

(2)  Let $(V, W, N_1,\dots, N_n, F)$ be an object of $DH_n$. Then $(V, W, \hat N_1, \dots, \hat N_n, \hat F)$, where $\hat N_j$ are as in \ref{Dthm2} and $\hat F$ is as in \ref{hatF1}, 
is an object of $\hat DH_n$.

\end{sbprop}

We call the object of $\hat D_{n,E}$ (resp. $\hat DH_n$) associated to an object of $D_{n,E}$ (resp. $DH_n$) in \ref{assl} , {\it the associated SL(2)-orbit}.

The proof of  \ref{assl} (1) is easy (the key point is (v) in \ref{Dthm2}).

The following  counterpart for $DH_n$ of  (v) in \ref{Dthm2} proves \ref{assl}  (2). 

\begin{sbprop}\label{hatF2} Let $(V, W, N_1,\dots, N_n, F)$ be an object of $DH_n$  with associated $\tau=(\tau_j)_{0\leq j\leq n}$.  Let   $\hat F$ be as in \ref{hatF1}.  Then we have 
$\tau_j(a)\hat F=\hat F$ for any $0\leq j\leq n$ and any $a$.

\end{sbprop}
In the case of IMHM, this is Lemma 5.5 in \cite{BP}.  We give the proof in the general case in \ref{pfhatF1} and \ref{pfhatF2}  below after preparations.

\begin{sbpara}\label{D2sp} Let $(V, W, N, F)$ be an IMHM of one variable.
Then $(V, W, \exp(iN)\hat F)$ is a mixed Hodge structure. Let $\tau'_0$ be the representation of $\mathbb{G}_m$ on $V$ defined by the canonical splitting of $W$ associated to the mixed Hodge structure. On the other hand, let  $(V, W, N, \alpha)$ be the  Deligne system of one variable associated to $(V, W, N, F)$ (Section 2.2) and consider its $\tau_0$.

(1) An important theorem of Deligne is 
$$\tau'_0=\tau_0.$$
This is introduced in \cite{BP} as Lemma 2.2 and the proof is given in that paper. 

(2) On the other hand, in \cite{KNU1}, it is proved that $\tau'_0(a)\hat F =\hat F$ for any $a$. 

By (1) and (2),  we have $\tau_0(a)\hat F=\hat F$. 

\end{sbpara}

\begin{sblem}\label{3.3lem1} Let $(V, W, N, F)$ be a DH system of one variable. Assume $W$ is pure and $\hat F=F$. Then this object is an SL(2)-orbit.

\end{sblem}

This is evident. 

The following is a special case of Theorem \ref{thm1} (this theorem shows that the assumption $\hat F=F$ is not necessary in the following lemma).
\begin{sblem}\label{3.3lem2}  Let $(V, W, N, F)$ be a DH system of one variable. Assume $\hat F=F$. Then this object is an IMHM.

\end{sblem}

This follows from Lemma \ref{3.3lem1}. 

\begin{sbpara}\label{pfhatF1}
We prove Proposition \ref{hatF2} in the case $n=1$. 
Let 
$(V, W, N, F)$ be a DH system of one variable. Then   $(V, W, N, \hat F)$ is a DH system of one variable  and satisfies the assumption of Lemma \ref{3.3lem2}. 
Hence it is an IMHM by Lemma \ref{3.3lem2}. Hence by \ref{D2sp}, we have $\tau_0(a)\hat F=\hat F$.

\end{sbpara}

\begin{sbpara}\label{pfhatF2} We prove Proposition \ref{hatF2} in general by induction on $n$. Assume $n\geq 2$.  Note that
$(W^{(1)}, N_2, \dots, N_n, F)$ is an object of $DH_{n-1}$ and the  associated action  $(\tau'_j)_{0\leq j\leq n-1}$ of $\mathbb{G}_m^n$ is given by $\tau'_j=\tau_{j+1}$. 
By the hypothesis of induction,  
 $(V, W^{(1)}, \hat N_2, \dots, \hat N_n, \hat F)$ is an SL(2)-orbit. From this, we have
 
 (1) $(V, W^{(1)}, F')$ with $F'=\exp(\sum_{j=2}^n i\hat N_j)\hat F)$ is a mixed Hodge structure. 
 
 (2) $\tau_1(a)F'=F'$ for any $a$. 
 
(1) and (2) hold also when we replace $V$ by $U:=W_w$ ($w\in \Z$) and replace $W^{(1)}$, $N_j$ and $F$ to their restrictions to $U$.  From this, we see that  $(V, W, N_1, F')$ is a DH system of one variable. By the case $n=1$ of 
\ref{hatF2} proved in \ref{pfhatF1},  and by (2) which shows that the functor $F\mapsto\hat F$ applied to $F'$ does not change $F'$, we have $\tau_0(a)F' =F'$ for any $a$. 
Since $\tau_0(a)\hat N_j\tau_0(a)^{-1}=\hat N_j$  for any $j$, this proves that $\tau_0(a)\hat F=\hat F$. This completes the proof of \ref{hatF2} and hence the proof of \ref{assl}. 
\end{sbpara}

\begin{sbprop}\label{convsl}  For an object of $D_{n,E}$ with $E=\R$ or $\C$, or for an object of $DH_n$, we use the notation
$$\beta(y)=\prod_{j=0}^n \tau_j((y_j/y_{j+1})^{1/2})\quad \text{for}\;\; y=(y_0,\dots, y_n)\in \R_{>0}^{n+1}$$
where $y_{n+1}$ denotes $1$.

(1) For an object of $D_{n,E}$ with $E=\R, \C$ or of $DH_n$,  for $y=(y_j)_{0\leq j\leq n}\in \R_{>0}^{n+1}$ such that $y_j/y_{j+1}\to \infty$ ($0\leq j<n$), 
we have the convergences  $$\beta(y)y_jN_j\beta(y)^{-1}\to \hat N_j, \quad  \beta(y)(\sum_{j\in I} y_jN_j)\beta(y)^{-1}\to \sum_{j\in I} \hat N_j$$ for any subset $I$ of $\{1,\dots,n\}$.

(2) For any object of $DH_n$, for $y=(y_j)_{0\leq j\leq n}\in \R_{>0}^{n+1}$ such that $y_j/y_{j+1}\to \infty$ ($0\leq j\leq n$),  we have the convergences  $$\beta(y)F\to \hat F, \quad 
\beta(y)\exp(\sum_{j\in I}  iy_jN_j)F\to \exp(\sum_{j\in I} i\hat N_j)\hat F$$ for any subset $I$ of $\{1,\dots,n\}$. 
\end{sbprop}
This is proved easily.

\begin{sbrem} The terminology  SL(2)-orbit in this paper is different from that in \cite{KNU}. In \cite{KNU}, we called an IMHM $(V, W, N_1, \dots, N_n, F)$ an SL(2)-orbit 
if $$\tau_k(a)N_j\tau_k(a)^{-1}=N_j\;\;\;(1\leq k<j\leq n)\quad \tau_k(a)F=F\;\;\;(1\leq k \leq n)$$
(for any $a$). (The difference is that  $1\leq k$ is in this  condition though $0\leq k$ is in the condition for SL(2)-orbit  in this paper.)  If $n=0$, this condition is empty. (In the case of SL(2)-orbit  of DH system in this paper, we require in the case $n=0$ that $F=\hat F$. )
If $n\geq 1$, this condition is equivalent to the condition 
$$\tau_k(a)N_j\tau_k(a)^{-1}=N_j\;\;\;(2\leq j\leq n, \;  0\leq k<j), \quad   \tau_k(a)F=F\;\;\;(0\leq k \leq n).$$
(The difference from the condition of this paper is that there is no condition on $N_1$ here.)

In \cite{KNU}, the associated SL(2)-orbit meant $(N_1, \hat N_2, \dots, \hat N_n, \hat F)$ ($N_1$ appears in place of $\hat N_1$ of this paper) if $n\geq 1$ and meant $F$ if $n=0$. 

In the pure case, there is no difference in these two formulations of SL(2)-orbits. 

The formulation  of SL(2)-orbit  in \cite{KNU}  is
useful for the study of classifying spaces of degenerating mixed Hodge structures (see \cite{KNU}). 

\end{sbrem}

\section{Proofs of theorems}

\subsection{A Deligne-Hodge system generates an IMHM}

We prove  Theorem \ref{thm1} in Introduction. We also prove
 
\begin{sbthm}\label{RMF} Let $E=\R$ or $\C$, and let $(V, W, N_1,\dots, N_n, \alpha)$ be an object of $D_{n,E}$. Take $1\leq \ell\leq  j<k\leq n$. Then for $y_t>0$ ($\ell+1\leq t\leq k$) such that $y_t/y_{t+1}\gg 0$ ($\ell+1\leq  t<k$), 
$W^{(k)}$ is the relative monodromy filtration of $\sum_{t=\ell+1}^k y_tN_t$ with respect to $W^{(j)}$. In other words, $(V, W^{(j)}, \sum_{t=\ell+1}^k y_tN_t, \tau_k)$ is a Deligne system of one variable. 

\end{sbthm}

\begin{pf} 

For $y=(y_t)_{0\leq t\leq n}\in \R_{>0}^{n+1}$, let  $N_y=\beta(y) (\sum_{t=\ell+1}^k  y_tN_t)\beta(y)^{-1}$ where $\beta(y)$ is as in \ref{convsl}. Let  $\hat N= \sum_{t=\ell+1}^k \hat N_t$. Then by \ref{convsl},  $N_y$ converges to $\hat N$ when $y_t/y_{t+1}\to \infty$ ($0\leq t <n$). By \ref{sl2nm},  the map $\hat N^m : \gr^{W^{(k)}}_{w+m} 
\gr^{W^{(j)}}
_w\to \gr^{W^{(k)}}_{w-m}\gr^{W^{(j)}}_w$ is an isomorphism for any $w\in \Z$ and any $m\geq 0$. It follows that if $y_t/y_{t+1}\gg 0$ ($0\leq t<n$), the map $N_y^m:  \gr^{W^{(k)}}_{w+m} \gr^{W^{(j)}}_w\to \gr^{W^{(k)}}_{w-m}\gr^{W^{(j)}}_w$ is an isomorphism and hence the map $(\sum_{t=\ell+1}^k y_tN_t)^m : \gr^{W^{(k)}}_{w+m} 
\gr^{W^{(j)}}
_w\to \gr^{W^{(k)}}_{w-m}\gr^{W^{(j)}}_w$ is an isomorphism.
  \end{pf}

\begin{sbpara}

We prove Theorem \ref{thm1}.

By Theorem \ref{RMF}, the condition (h) in the definition of the notion IMHM is satisfied. In fact, for $N'_j= \sum_{k=0}^j a_{j,k} N_k$ with $a_{j,k}>0$ such that
 $a_{j,k}/a_{j, k+1}\gg 0$ ($1\leq k <j$), by Theorem \ref{RMF}, 
$W^{(j)}$ is the relative monodromy filtration of $N'_j$  with respect to $W$. 

It remains to consider the condition (g), i.e. the polarizability of $\gr^W$.
On gr$^W_w$, put the bilinear form in \ref{sl2pol}. 
For $y=(y_j)_{0\leq j\leq n}\in \R_{>0}^{n+1}$, let $F(y)=\exp(\sum_{j=1}^n iy_nN_n)F$, $I= \exp(\sum_{j=1}^n i\hat N_j)\hat F)$. 
Let  $\beta(y)$ be as in \ref{convsl}. Then by \ref{convsl},  $\beta(y)F(y)$ converges to $I$
when $y_j/y_{j+1}\to \infty $ ($0\leq j\leq n$, $y_{n+1}$ means $1$). 
Since $(V, W, I)$ is a mixed Hodge structure (\ref{sl2pol}), we have that $(V, W, \beta(y)F(y))$ is a mixed Hodge structure  when $y_j/y_{j+1}\gg 0$. 
Hence we can consider the Hermitian form associated to $(\langle\;,\;\rangle_w, \beta(y)F(y)(\gr^W_w))$. This Hermitian form converges to the Hermitian form associated to $( \langle\;,\;\rangle_w, I(\gr^W_w))$ which is positive definite.  
Hence the former Hermitian form is positive definite if $y_j/y_{j+1}\gg 0$. Hence when $y_j/y_{j+1}\gg 0$, $(V, W, F(y))$ is a mixed Hodge structure and $(\gr^W_w, \langle\;,\;\rangle_w, F(y)(\gr^W_w))$ is a polarized Hodge structure. 
This proves Theorem \ref{thm1}. 

\end{sbpara}

By Theorem \ref{thm1}, 
SL(2)-orbit theorems for IMHM in \cite{Scw}, \cite{CKS},  \cite{KNU1} are generalized to $DH_n$. For example (\cite{CKS}, Theorem 4.20 (vii)), we have

\begin{sbthm}  Let $(V, W, N_1,\dots, N_n, F)$ be an object of $DH_n$. Assume $W$ is pure. Then there is a convergent series $g(T_1,\dots, T_n)\in \End_{\R}(V)[[T_1,\dots, T_n]]$ with constant term $1$ such that 
when $y_j/y_{j+1}\gg 0$ ($y_{n+1}$ denotes $1$), we have
$$\exp(\sum_{j=1}^n iy_jN_j)F=  g(y_2/y_1, y_3/y_2, \dots, y_{n+1}/y_n)\cdot \prod_{j=1}^n \tau_j((y_{j+1}/y_j)^{1/2}))\cdot \exp(\sum_{j=1}^n i\hat N_j)\hat F   $$ 

\end{sbthm}

\subsection{On Theorem \ref{Dthm} and Theorem \ref{thm3}}

Concerning Theorem \ref{Dthm},
we give a more precise statement about the convergence of the splitting of $W$. 
\begin{sbthm}\label{Wspas}

Let $n\geq 1$ and let $(V, N_1, \dots, N_n, F)$ (resp. $(V, N_1, \dots, N_n, \alpha)$) be  an object of $DH_n$ (resp. $D_{n,E}$ with $E=\R$ or $\C$).  
For $y=(y_1,\dots, y_n)\in \R_{>0}^n$, let $H(y)=(V, W, \exp(\sum_{j=1}^n iy_jN_j)F)$ (resp. $H(y)=(V, W, \sum_{j=1}^n y_jN_j,\alpha)$). Let $\hat H=(V, W,  \exp(\sum_{j=1}^n i\hat N_n)\hat F)$ (resp. $\hat H=(V, W, \sum_{j=1}^n \hat N_j, \alpha)$). Let $h=n$ (resp. $h=n-1$). Then there is a family $(u_m)_{m\in\N^h}$ of 
 $\R$ (resp. $E$)-linear maps $u_m:V\to V$
having the following properties (i)--(iii).  

\medskip
(i) $u_0=1$. $u_mW_w\subset W_{w-1}$ for any $m\neq 0$ and any $w\in \Z$. For $1\leq j\leq n$, $u_mW^{(j)}_w\subset W^{(j)}_{w+m(j)}$ for any $m$ and any $w$. 

\medskip

(ii) Let $u(T_1, \dots, T_h)= \sum_{m\in \N^h} \;u_mT_1^{m(1)}\dots T_n^{m(h)}$. Then there is $c>0$ such that $u(T)$ absolutely converges if $|T_j|<c$ for all $j$.

\medskip
(iii) 
For $y_j>0$ ($1\leq j\leq n$) such that  $y_j/y_{j+1}\gg 0$ ($1\leq j\leq n$ where $y_{n+1}$ denotes $1$) (resp.  ($1\leq j<n$)), let 
$s(y):\gr^W\overset{\cong}\to V$ be the splitting of $W$ associated to the mixed Hodge structure (resp. Deligne system of one variable) $H(y)$. Let $\hat s: \gr^W\overset{\cong}\to V$ be the splitting of $W$  associated to the mixed Hodge structure (resp. Deligne system of one variable) $\hat H$. Then there is $c>0$ such that
$$s(y)= u(y_2/y_1, \dots, y_{h+1}/y_h)\hat s$$
when $y_j/y_{j+1}\gg 0$ ($1\leq j\leq h$).

\end{sbthm}

By Theorem \ref{thm1}, Theorem \ref{Dthm} and Theorem \ref{Wspas} for $DH_n$  follow from the corresponding result \cite{KNU1}, Theorem 0.5 for IMHM.

We will prove  in \ref{pfsp} the part concerning $D_{n,E}$ ($E=\R, \C$)
 by reducing it to the part on $DH_n$.

\begin{sblem}\label{spsp}  Let $E=\R$, let $(V, W, N, \alpha)$ be a Deligne system of one variable with the associated $(\tau_j)_{j=0,1}$ and let  $(V^{\oplus 2}, W^{\oplus 2}, N^{\oplus 2}, F)$ be the corresponding object of $DH_1$ (Section 2.3). Then 
$(V, W^{\oplus 2}, \exp(iN^{\oplus 2})F)$ is a mixed Hodge structure, and we have $\tau'_0=\tau_0^{\oplus 2}$ where $\tau'_0$ denotes the splitting 
of $W^{\oplus 2}$ associated to this mixed Hodge structure.  

\end{sblem}

\begin{pf} By Theorem \ref{thm1}, any object of $DH_1$ is an IMHM. Hence $(V^{\oplus 2}, W^{\oplus 2}, N^{\oplus 2}, F)$ is an IMHM. It  is easy to see $\hat  F=F$. Hence the result  follows from  \ref{D2sp} (1).
\end{pf}

\begin{sbpara}\label{pfsp}  We prove the parts of $D_{n,E}$ in Theorem \ref{Dthm} and Theorem \ref{Wspas}.

By Lemma \ref{spsp}, these theorems for $D_{n,E}$ are reduced to the parts  for $DH_n$ by using the functor $D_{n,E}\to DH_n$ (Section 3.3). In fact, we have the result $s(y)^{\oplus 2}=u(y)\hat s^{\oplus 2}$ with $u(y)\in \End_{\R}(V^{\oplus 2})$. If we write $u(y)$ in the form $(u(y)_{s,t})$ ($s,t\in \{1,2\}$) of $(2,2)$-matrix in which all entries are elements of $\End_{\R}(V)$, $s(y)^{\oplus 2}
=u(y)\hat s^{\oplus 2}$ shows that 
$s(y)= (u(y)_{1,1}+u(y)_{1,2})\hat s$, and hence we have the result $s(y)=u(y)\hat s$ with $u(y)\in \End_{\R}(V)$. If $E=\C$, by replacing $u$ by $(1/2) (u \circ i + i \circ u)$ where $i$ denotes the multiplication by $i$, we can take $u(y)\in \End_{\C}(V)$. 
Finally, since $s(y)$ depends only on the ratio $(y_1:\dots:y_n)$, the series $u(y)$ in $y_2/y_1, \dots, y_{n+1}/y_n$ ($y_{n+1}$ denotes $1$) is actually series in $y_1/y_2, \dots, y_n/y_{n-1}$. 

\end{sbpara}

\begin{sbpara}\label{atwi}

For $a\in \R$, define the functor $\theta^a:DH_n\to DH_n$ as $$(V, N_1, \dots, N_n, F)\mapsto (V, N_1', \dots, N'_n, F) \quad\text{where}\;\;N'_j = \sum_{k=0}^{j-1} (a^k/k!)N_{j-k}.$$
That is, $(N'_1, \dots, N'_n)^t=\exp(aR) (N_1, \dots, N_n)^t$ where $R$ is the $(n,n)$ matrix whose $(j,k)$-th entry is $1$ if $k=j-1$ and  $0$ otherwise. 

For an object $H$ of $DH_n$, we have 
\smallskip

(1) $\theta^{a+b}H=\theta^a(\theta^bH)$.
\smallskip

Since $\theta^a\theta^{-a}$ is the identity functor by (1), we see that
\smallskip

(2) $\theta^a:DH_n\to DH_n$ is an equivalence of categories.
\smallskip

By Theorem \ref{thm1}, we have
\smallskip

(3) If $H$ is an object of $DH_n$,  $\theta^aH$ is an  IMHM if $a\gg 0$. 
 \smallskip

For $a\in \R$, let  $DH_n^{(a)}$ be the full subcategory of $DH_n$ consisting of all  objects $H$ such that $\theta^aH$ is an IMHM. By (3), we have
\smallskip

(4) $DH_n = \cup_a\;  DH_n^{(a)}$. Note that $DH_n^{(a)}\subset DH_n^{(b)}$ if $a\leq b$.

\end{sbpara}

\begin{sbpara} We prove Theorem \ref{thm3}. 

First we prove the part concerning $DH_n$. Theorem \ref{thm3} is true if $DH_n$ is replaced by the category of IMHM of $n$ variables (\cite{Kas}).
For $a\in \R$, the category $DH_n^{(a)}$ in \ref{atwi}  is equivalent to the category of IMHM of $n$ variables by the functor $\theta^a$. This shows that $DH_n^{(a)}$ is an abelian category and the kernel and the cokernel are described as in Theorem \ref{thm3}. By (4) in \ref{atwi}, this proves Theorem \ref{thm3} for $DH_n$.

We prove the part concerning $D_{n,E}$. 
First we show that we can assume $E=\C$. This is because an object of $D_{n,E}$ or a morphism of $D_{n,E}$ comes from $D_{n,E'}$ for some subfield $E'$ of $E$ which is finitely generated over $\Q$. Then we have an embedding of $E'$ into $\C$ as a subfield. Hence by (1) of  Lemma \ref{EandE'}, we are reduced to the case $E=\C$.

Next by (2) of Lemma \ref{EandE'},  we can assume $E=\R$. 

We prove Theorem \ref{thm3} in the case $E=\R$.  We denote the functor $DH_n\to D_{n,\R}$ in Section 3.2 by $a$ and the functor $D_n\to DH_n$ in Section 3.3 by $b$. 
Let $f: A=(V_A, W_A,N_{1,A},\dots, N_{n, A}, \alpha_A)\to B=(V_B, W_B, N_{1,B},\dots, N_{n, B},\alpha_B)$ be a morphism of $D_{n,\R}$, let $V_K$ (resp. $V_C$) be the kernel (resp. cokernel) of $f: V_A\to V_B$, and let $W_K, N_{j, K}, \alpha_K$ on $V_K$  (resp. $W_C, N_{j,C}, \alpha_C$ on $V_C$) be the ones induced from those of $A$ (resp. $B$). Then $f$ induces a morphism $b(f): b(A) \to b(B)$ of $DH_n$, and the kernel and the cokernel of $b(f)$ are described as in the part of Theorem  \ref{thm3}   concerning $DH_n$ which we have proved. By applying the functor $a$, we see that $(V_K^{\oplus 2}, W_K^{\oplus 2}, N_{1,K}^{\oplus 2}, \dots, N_{n, K}^{\oplus 2}, \alpha_K^{\oplus 2})$ (resp. $(V_C^{\oplus 2}, W_C^{\oplus 2}, N_{1,C}^{\oplus 2}, \dots, N_{n, C}^{\oplus 2}, \alpha_C^{\oplus 2})$) is an object of $D_{n,\R}$. This shows that 
$K:=(V_K, W_K, N_{1,K}, \dots, N_{n, K}, \alpha_K)$ (resp. $C:=(V_C, W_C, N_{1,C}, \dots, N_{n, C}, \alpha)$) is an object of $D_{n,\R}$. We have shown that the kernel and the cokernel of a morphism in $D_{n,\R}$ exist and are described as in Theorem \ref{thm3}. Let $I$ be the cokernel of $K\to A$ (resp. $J$ be the kernel of $B\to C$). Since $DH_n$ is an abelian category, the canonical morphism  $b(I)\to b(J)$ is an isomorphism. By applying the functor $a$, we see that the canonical morphism $I^{\oplus 2}\to J^{\oplus 2}$ is an isomorphism. Hence the canonical morphism $I\to J$ is an isomorphism. This proves Theorem \ref{thm3} for $D_{n,\R}$. 

\end{sbpara}

Kazuya KATO,

Department of mathematics,
University of Chicago,
Chicago, Illinois, 60637,
USA,

\medskip

Running title. SL(2)-orbit theorems


\begin{thebibliography}{99}



\bibitem{BK}{\sc Bloch, S, Kato, K.},
\newblock
{\em Asymptotic behaviors of heights and regulators in degeneration}
in preparation.




\bibitem{CKS}{\sc Cattani, E., Kaplan, A., Schmid, W.},
\newblock
{\em Degeneration of Hodge structures},
Ann. of Math. {\bf 123} (1986), 457--535.





\bibitem{De}
{\sc Deligne  ~P.},
{\em La conjecture de Weil},
Publ. Math. I.H.E.S. {\bf 52}
(1980), 137--252.



\bibitem{Kas}{\sc Kashiwara, M,},
\newblock
{\em A study of variation of mixed Hodge structure}, 
Publ. RIMS, Kyoto Univ. {\bf 22},  (1986),  991--1024.


\bibitem{KK}{\sc Kato, K.},
\newblock
{\em    p-Adic period domains and toroidal partial compactifications, I},
Kyoto J. Math. {\bf 51}, (2011), 561-631.

\bibitem{KNU1}{\sc Kato, K., Nakayama, C., Usui, S.},
\newblock
{\em SL(2)-orbit theorem for degeneration of mixed Hodge structure}

\bibitem{KNU},
{\sc Kato, K., Nakayama, C., Usui, S.},
{\em Classifying spaces of degenerating mixed Hodge structures},
I., Adv. Stud. Pure Math. {\bf 54}, 187--222, II., Kyoto J. Math. {\bf 51}, 
149--261, III., J. of Algebraic Geometry, {\bf 22}, (2013), 671--772.


\bibitem{BP} 
{\sc Brosnan, P, Pearlstein, G.}, 
{\em On the algebraicity of the zero locus of an admissible normal function},
to appear in Composition Math.  arXiv:0910.0628 

\bibitem{Scw}{\sc Schmid, W.},
\newblock
{\em  Variation of Hodge structure:
the singularities of the period mapping}
 Invent. Math.
{\bf 22} (1973), 
211--319.


\bibitem{Scc}
{\sc Schwarz, C.}, 
{\em Relative monodromy weight filtrations},
 Math. Z. {\bf 236} (2001), 11--21.
 
 
\end{thebibliography}
\end{document}